\magnification=1100
\baselineskip=16pt

\hsize=152mm
\vsize=205mm

\let \al=\alpha
\let \be=\beta
\let \var=\varphi
\let \vare=\varepsilon

\let \de=\delta

\let \la=\lambda

\let \ga=\gamma
\let \k=\kappa
\let \p=\partial

\let \med=\medskip
\let \smal=\smallskip
\let \dps =\displaystyle
\let \q=\quad

\let \ol =\overline
\let \ul=\underline

\def\R{{\rm I\kern
-1.6pt{\rm R}}}
\def\C{{\rm |\kern
-4.6pt{\rm C}}}
\def\N{{\rm I\kern
-4.0pt{\rm N}}}

\def\eq{\eqno}
\def\q{\quad}

\def\ter{\hfill \vrule width 5 pt height 7 pt depth - 2 pt\smallskip}

\def\system#1{\left\{\null\,\vcenter{
\ialign{\strut\hfil$##$&$##$\hfil&&\enspace$##$\enspace&
\hfil$##$&$##$\hfil\crcr#1\crcr}}\right.}

\centerline{\bf Persistence and permanence for a class of  functional differential equations} 
\centerline{\bf with infinite delay} 

\vskip0.3in

\centerline{Teresa Faria\parindent=0cm{\footnote {$^{a}$}{Departamento de Matem\'atica  and CMAF, Faculdade de
Ci\^encias, Universidade de Lisboa, Campo Grande, 1749-016 Lisboa,
Portugal}}
\parindent=0cm{\footnote{$^{\star}$}{
Fax:~+351 21 795 4288, tel:~+351
21 790 4929, e-mail:~teresa.faria@fc.ul.pt.}}
}

\vskip0.3in

{\it \centerline{Dedicated to Professor John Mallet-Paret,}
\centerline{ on the occasion of his 60th birthday}}

\vskip0.3in

{\it Suggested running head: Persistence and permanence for FDEs with infinite delay}

\vskip0.3in

\centerline{\bf Abstract}
\vskip0.1in

The paper deals with a class of cooperative functional differential equations  (FDEs) with infinite delay, for which sufficient conditions for persistence and permanence are established. Here, the persistence refers to all  solutions with initial conditions that are positive, continuous and bounded. 
The present method
applies to a very broad class of abstract systems of FDEs with infinite delay, both autonomous and non-autonomous, which include many important models used in mathematical biology.  Moreover,   the  hypotheses imposed are in general very easy to check. The results are illustrated with some selected examples.

\vskip0.15in

{\parindent=0cm{\it Keywords}:  infinite delay,  persistence, permanence, quasimonotone condition,  Lotka-Volterra model.

\smal

}

\

{\bf 1. Introduction}

\med



This paper is concerned with the persistence and permanence for a class of functional differential equations (FDEs) with infinite delay written in abstract form as
$$
x_i'(t)=F_i(t,x_t)-x_i(t)G_i(t,x_t), \quad i=1,\dots,n
\eq(1.1)$$
where    $F=(F_1,\dots,F_n),G=(G_1,\dots,G_n):D \to \R^n$ are  continuous on  $D\subset \R\times {\cal B}$, and ${\cal B}$ is an adequate Banach space of continuous functions defined on $(-\infty,0]$ with values in $\R^n$. As usual, $x_t$ denotes the entire past history of the system up to time $t$, or, in other words, $x_t(s)=x(t+s)$ for $s\le 0$.  Special emphasis will be given to the {\it autonomous} version of (1.1),
$$
x_i'(t)=F_i(x_t)-x_i(t)G_i(x_t), \quad i=1,\dots,n
\eq(1.2)$$
where    $F_i,G_i:\Omega\subset{\cal B}\to\R$ are  continuous for all $i$. 

We shall consider (1.1)  as a general framework for some models from mathematical  biology, therefore  only  positive (or non-negative) solutions of (1.1) are meaningful. Initial conditions at a fixed time  $\sigma$ will be taken to be $x(s)=\var(s),\, s\le \sigma$, with  $\var:(-\infty,0]\to \R^n$ continuous, bounded and  positive  (or non-negative). 


The family of  FDEs (1.1) is required to satisfy the so-called quasimonotone condition [30] -- here, we abuse the standard terminology, and  these systems will  be also called {\it cooperative systems} (see [30] for a rigorous definition).  
The main technique employed here is the theory of monotone dynamical systems for retarded FDEs well-established in [29,30,34]: the quasimonotone condition  implies that the solution operator is monotone relative to the initial condition function $\var$ (but not necessarily strongly monotone) and allows comparison of solutions between related FDEs.  However, it is not always clear how to apply the tools of this theory to the situation of infinite delays. In particular, for (1.1) it is not obvious which conditions should be imposed in order to guarantee the persistence of all  solutions with positive, bounded initial conditions. 

From earlier works (see e.g. [16,25,26]), it has been apparent  that a  rigorous treatment of FDEs with infinite delay requires a particular attention to its abstract formulation.
In Section 2, we choose an {\it admissible} phase space  ${\cal B}$ to deal with (1.1), in the sense that it should satisfy some fundamental axioms introduced by Hale and Kato [17].  To address the question of permanence,  we shall  need  bounded positive orbits to have compact closure in ${\cal B}$, so, in addition,  we choose a state space that is a {\it fading memory}  space [14,18,27].  
The contents and organization of the remainder of paper are  described below. 

Section 3 is the core of the paper, where the main theorems  are established, referring to  persistence (Theorem 3.1) and permanence (Theorem 3.2) of autonomous cooperative systems (1.2). A key point of our research is that quite reasonable sufficient conditions are enough to guarantee the persistence (or even the uniform persistence) of system (1.2) in the entire set of solutions with  initial conditions that are bounded, continuous and positive. In Section 4, we use the results and techniques in Section 3 to further deduce some persistence and permanence criteria for non-autonomous cooperative systems (1.1). We emphasize that, both in Section 3 and 4,  no global asymptotical constancy is invoked to establish the persistence or permanence of a system. 


Section 5 is dedicated to  applications and  is divided into three subsections, each  dealing with a selected  important  example from population dynamics. 
In the first example,   sufficient conditions are established for the persistence and permanence of a cooperative Lotka-Volterra model with infinite delay and patch structure, improving significantly some results  in a former paper of the author [10].  In fact, the study of this Lotka-Volterra model was a main source of motivation for the present study, since there was an open problem concerning the persistence left unsolved in  [10]. For alternative techniques for competitive Lotka-Volterra system, see e.g. [1,8,24,32,33].
%
The second example refers to an FDE modelling the growth of a single or multiple species divided into $n$ classes and  following a modified delayed logistic law; the system has unbounded discrete time-varying delays, also includes   dispersal terms among the classes,  and may be interpreted as a generalization of the modified scalar delayed logistic equation proposed  by Arino et al.  [6]. The last selected example concerns an FDE   system for the interaction of two species  structured into mature and immature classes, with the two adult populations in competition; this system  is based on Aiello and Freedman's  model for a single species [2], and the finite delay case was studied by  Al-Omari and Gourley [5], among others. Although the system is competitive,   the method in [5] relies on the construction of  sequences of auxiliary cooperative systems, which are used to prove the global attractivity of a positive equilibrium.
   Here, we show that Al-Omari and Gourley's method  can be extended to the case of infinite distributed delay. 
These examples  have been considered by many authors, and additional references  will be given in Section 5.
Many other examples from the literature could be analyzed, but the main purpose of Section 5 is  to illustrate the application of our main results.

\bigskip

{\bf 2. Preliminaries: phase space and notations}

\smal

In this preliminary section, we set an abstract framework to deal with  FDEs with infinite delay.  In view of the  unbounded delays, the phase space ${\cal B}$ should satisfy some fundamental axioms which guarantee that  the classical results of existence, uniqueness, continuation, and continuous dependence  of solutions  are valid -- a subject well establish in the literature. Secondly, for our purposes, it is desirable that   positive orbits of bounded solutions (with bounded derivatives) are precompact in ${\cal B}$.  A convenient choice of ${\cal B}$ is set below.

Let $g$ be a
function satisfying the following properties:
\smal

 {(g1)} $g:(-\infty ,0]\to [1,\infty)$ is a nonincreasing continuous function, $g(0)=1$;

 {(g2)} $\displaystyle{\lim_{u\to 0^-}{{g(s+u)}\over {g(s)}}=1}$ uniformly on $(-\infty ,0]$;

 {(g3)} $g(s)\to \infty$ as $s\to -\infty$.

\noindent 
For $n\in\N$, define the Banach space ${\cal B}=C_g^0$,
$$C_g^0=C_g^0(\R^n):=\left\{ \phi\in C((-\infty, 0];\R^{n}) : \lim_{s\to -\infty }{{|\phi(s)|}\over {g(s)}}=0\right\},$$
 with the norm
$$\|\phi\|_g=\sup_{s\le 0}{{|\phi(s)|}\over {g(s)}},$$
where $|\cdot|$ is any chosen norm  in $\R^n$. 
To fix terminology we suppose that $\R^n$ is equipped with the supremum norm, $|x|=|x|_\infty=\max\limits_{1\le i\le n}|x_i|$, for $x=(x_1,\dots,x_n)\in\R^n$. 
Consider also  the space $BC=BC(\R^n)$ of bounded  continuous functions $\phi:(-\infty, 0]\to
\R^n$.
Here, $BC$ is considered as
a subspace of some space $C_g^0$, so $BC$ is endowed with the norm of $C_g^0$. The more explicit notation $BC_g$ will  be  also used.

The space $C_g^0$ is an admissible phase space [17,18] for $n$-dimensional FDEs with infinite delay
  written in the abstract form
$$
x'(t)=f(t,x_t),
\eq(2.1)$$
where $f:D\subset \R\times C_g^0\to \R^n$ is continuous and, 
as usual, segments of solutions in the phase space $C_g^0$ are denoted by $x_t$: $x_t(s)=x(t+s), s\le 0$. When $f$ is regular enough,  
the initial value problem is well-posed, in the sense that for each $(\sigma,\var)\in D$ there exists a unique solution $x(t)$ of  the problem $x'(t)=f(t,x_t), x_\sigma=\var$,  defined on a maximal interval of existence. This solution will be denoted by $x(t;\sigma,\var)$ in $\R^n$ or  $x_t(\sigma,\var)$ in $C_g^0$. For autonomous systems $x'(t)=f(x_t)$ and $\var\in C_g^0$, the solution of $x'(t)=f(x_t),\, x_0=\var$ is often simply denoted by $x(t;\var)\in\R^n$ and $x_t(\var)\in C_g^0$.  When considering more than one FDE of the form $x'(t)=f(t,x_t)$ or $x'(t)=f(x_t)$, the notations $x(t;\sigma,\var,f)$ or $x(t;\var,f)$ will be also used, where the argument $f$ is included to avoid any confusion over which FDE is being considered. The space $C_g^0$ has important properties, namely it is a {\it fading memory space} -- although  not always a {\it uniform}  fading memory space. For definitions and results, see [14,18,27].  For FDEs $x'(t)=f(x_t)$ in fading memory spaces,   positive orbits of  solutions $x(t)=x(t;\var)$    with {\it bounded} initial conditions (i.e., for $\var \in BC$)  which are bounded and have bounded derivatives on $[0,\infty)$ are {\it precompact} in $C_g^0$ -- a property that will be used to prove our results of permanence. 

 Alternatively, one may consider the phase space  ${\cal B}=UC_g$, where
$UC_g=UC_g(\R^n):=\big\{ \phi\in C((-\infty, 0];\R^{n}) : \sup_{s\le 0}{{|\phi(s)|}\over {g(s)}}<\infty,
{{\phi(s)}\over {g(s)}}\ {\rm is\ uniformly\ continuous\ on}\ (-\infty, 0]\big\},$
with the norm $\|\cdot\|_g$ defined above. The space $UC_g$ is an admissible phase space, in the sense of Hale and Kato [17], but it is not necessarily a fading memory space -- contrary to what is stated in [19, p.~47]. In fact, Haddock and Hornor [14] completely described the functions $g$ satisfying (g1)-(g2) for which $UC_g$ is a fading memory space; in particular, $g$ must satisfy the condition $e^{-\gamma_1 s}\le g(s)\le e^{-\gamma_2 s}$ for $s\le -M$, for some $\gamma_1,\gamma_2,M>0$.  For example, (g1)-(g3) hold for $g(s)=1-s$, but $g$ does not satisfy this latter condition. If $g(s)=e^{-\gamma s}, s\le 0,$ for some $\gamma >0$, the spaces $C_g^0$ and $UC_g$ are uniform fading memory spaces, and are often denoted by $C_\ga^0$ and $UC_\ga$, respectively.
Some of the results established here require having precompactness of bounded positive orbits, therefore there is a clear  advantage in considering the  phase space $C_g^0$, rather than the usual $UC_g$.

We now set some notation. A vector $c$  in $\R^n$ is said to be {\it positive} if all its components are positive, and we write $c>0$. A function $\var:(-\infty,0]\to \R^n$ is said to be {\it positive}, with notation $\var>0$,    if the vectors $\var(s)$ are positive for all $s\le 0$. We define and denote in a similar way    {\it non-negative} vectors and {\it non-negative} functions.  In the space $C_g^0$ (or $UC_\ga$), a vector $c$  is identified with the constant function $\psi(s)=c$ for $s\le 0$. The non-negative cones  of $\R^n$ and $BC$ are  $\R^n_+=\{x\in\R^n:x\ge 0\}$ and $BC_+=BC_+(\R^{n})=\{ \var\in BC: \var(s)\ge 0$ for all $s\le 0\}$, respectively.  Motivated by the applications to mathematical biology, we shall take the set
$$BC_0=\{ \var\in BC: \var(s)>0\ {\rm for \ all\ } s\le 0\}$$ 
as the set of admissible initial conditions for (2.1). However, in applications, frequently  more general orders on  ${\cal B}$ are considered,  induced by  other cones $K\subset {\cal B}$. 
  We use   $(x)_i$ or simply $x_i$ to designate the component $i$ of a vector $x\in\R^n$. If $f$ is a function with values in $\R^n$,  it is understood that $f_i$ means the $i$th-component of $f$, for $1\le i\le n$.



\bigskip

{\bf 3. Main results}

\med

In this section, we prove the main results of the paper, about persistence and permanence for a class of autonomous FDEs of the form
$$
x_i'(t)=F_i(x_t)-x_i(t)G_i(x_t), \quad i=1,\dots,n,
\eq(3.1)$$
where    $F=(F_1,\dots,F_n),G=(G_1,\dots, G_n):BC\to \R^n$ are  continuous.  It is understood that $BC=BC_g$ is taken as a subspace of ${\cal B}=C_g^0$ (or ${\cal B}=UC_g$).  As mentioned above, we are concerned with the solutions $x(t)=x(t;\var)$ of (3.1) with initial conditions $$x_0=\var\q {\rm with}\q \var\in BC_0.
\eq(3.2)$$

We start with a preliminary lemma.

\proclaim{Lemma 3.1}. Consider $F,G:BC\to \R^n$ continuous. Then:\vskip 0cm
(i) if $F,G$ are bounded on  bounded sets of $BC$, and $x(t)$ is a non-continuable solution of (3.1)  defined on $[0,a)$ with $a<\infty$, then $\limsup_{t\to a^-}|x(t)|=\infty$;\vskip 0cm
\vskip 0cm
(ii) if  $F_i(\var)\ge 0$ for all $\var \in BC_+$ with $\var_i(0)=0\, (1\le i\le n)$, then the solutions   $x(t)$ with initial conditions  $\var\in BC_+$ satisfy $x(t)\ge 0$ for $t\ge 0$, whenever they are defined.

{\it Proof}. The statement in (i) is a classical result [17,18]; the proof of (ii) follows from [30,34].\ter


\med

System (3.1) reads as $$x'(t)=f(x_t),
\eqno (3.3)$$
where 
$$f=(f_1,\dots,f_n),\ f_i(\phi)=F_i(\phi)-\phi_i(0)G_i(\phi)\q {\rm for} \q \phi\in BC,1\le i\le n.$$
 
In the remainder of this section, $F=(F_1,\dots,F_n),  G=(G_1,\dots, G_n)$ are assumed to be continuous, bounded on bounded sets of $BC$, and regular enough so that the initial value problems (3.1)-(3.2)  have unique solutions -- which is the case if $f$ is Lipschitz continuous on bounded sets.

The main assumptions that will be  imposed, either in part or in total, are described below.

\vskip 3mm

{\parindent=1.5cm{

\item{(A1)} 
 for $\phi,\psi\in BC_+,\phi\le \psi$ and $\phi_i(0)=\psi_i(0)$, then (i) $F_i(\phi)\le F_i(\psi)$ and (ii) $G_i(\phi)\ge G_i(\psi), \ i=1,\dots,n$;

  \item{(A2)} for each $\de>0$ and each $ i=1,\dots,n$, there exists  $\vare>0$ such that for all $\phi,\psi\in BC_+,\phi\le \psi$ with $\|\psi-\phi\|_g<\vare$ and $\phi_i(0)=\psi_i(0)$, then $G_i(\phi)-G_i(\psi)<\de$;

\item{(A3)} there exists a constant $c>0$ such that for $\phi,\psi\in BC_+,\phi\le \psi$ and $\phi_i(0)=\psi_i(0)$, then $F_i(\phi)-F_i(\psi)\ge -c\|\psi-\phi\|_g, \ i=1,\dots,n$;

 \item{(A4)} $F$ is sublinear in $\R^n_+$, i.e., for $x\in \R^n_+$ and $\al\in (0,1)$, $F(\al x)\ge \al F(x);$

\item{(A5)} there exists a  vector $v\in\R^n_+$ such that $F(v)-Bv>0$, where 
$B=diag\, (G_1(0),\dots,G_n(0));$
\item{(A6)} there exists a  vector $q=(q_1,\dots,q_n)\in\R^n_+$  such that $F_i(q)-q_iG_i(Lq)<0$ for $L\ge 1, i=1,\dots,n$.

}}

\vskip 3mm

Some comments about the choice of the above hypotheses are in order. 

It is apparent that if the pair $(F,G)$ satisfies  one or more of the above hypotheses, the same happens for the pair $(\tilde F,\tilde G)$, where $\tilde F_i(\phi)=F_i(\phi)+G_i(0)\phi_i(0), \tilde G_i(\phi)=G_i(\phi)-G_i(0), 1\le i\le n$. Hence, without loss of generality, we can take $G(0)=0$, in which case   $B=0$ for the matrix $B$ in (A5). 
Assumption (A1) asserts that $F$ and $-G$, and thus $f$ in (3.3) as well, satisfy Smith's quasimonotone condition  given by (see [30])

\vskip 3mm

{(Q)} for $\phi,\psi\in BC_+,\phi\le \psi$ and $\phi_i(0)=\psi_i(0)$, then $f_i(\phi)\le f_i(\psi), \ i=1,\dots,n$.

\vskip 3mm

\noindent This condition implies that the semiflow for (3.3) is monotone on $BC_+$. Here we abuse the terminology, and also refer to systems satisfying (Q) as {\it cooperative systems}. 
Instead of (A2), we could simply demand $G$ to be uniformly continuous on $BC_+$, but this requirement is too strong for our purposes. Another advantage of the above formulation is that (A2) is trivially satisfied if each $G_i(\phi)$ depends only on $\phi_i(0)$.
It is clear that (A3) and (A4) are satisfied if $F$ is linear bounded; also, if $F$ is Lipschitz continuous on $BC_+$, then (A3) holds.
Condition (A6) is true whenever either $G\equiv 0$ and $F(q)<0$, or $\lim\limits_{L\to\infty}G_i(Lq)=\infty,\, i=1,\dots,n$, for some vector $q>0$. A careful analysis of the proofs below shows that, as alternatives to the sets of constraints (A4),(A5) and (A4),(A6), one can simply impose the hypotheses

\vskip 3mm

{\parindent=1.5cm{

\item{(A5')} there is  $v=(v_1,\dots,v_n)\in\R^n_+$ such that $F_i(lv)-lv_iG_i(lv)>0$ for $0<l\le 1,1\le i\le n,$

\vskip 2mm
\noindent and
\vskip 2mm

\item{(A6')} there is  $q=(q_1,\dots,q_n)\in\R^n_+$ such that $F_i(Lq)-Lq_iG_i(Lq)<0$ for $L\ge 1,1\le i\le n$,

}}

\vskip 3mm

\noindent respectively, without the additional requirement of the sublinearity of $F\big |_{\R^n_+}$. Although not essential,    the sublinearity condition (A4) allied to the quasimonotone condition (A1) implies, however, that solutions with non-negative  initial conditions remain non-negative. In fact, 
for $\var\in BC_+$ with $\var_i(0)=0$, from (A1) we have $f_i(\var)=F_i(\var)\ge F_i(0)$, $1\le i\le n$. On the other hand,  (A4) implies $F(0)\ge 0$. From Lemma 3.1(ii), it follows that $x(t;\var)\ge 0$ for $t\ge 0$ in its interval of definition.

 \smal

  For the definitions of persistence and permanence given below, see e.g. [19,31].
  
  \med

{\bf Definition 3.1}. A system $x'(t)=f(x_t)$ with $S\subset BC_+$ as set of admissible initial conditions is said to be {\it persistent}  (in $S$) if any  solution $x(t;\var)$ with initial condition $\var\in S$ is defined and bounded below from zero on $[0,\infty)$, i.e.,
$$\liminf_{t\to\infty} x_i(t;\var)>0,\q 1\le i\le n,$$ 
for any $\var\in S$. The system is said to be {\it permanent} (in $S$)  if it dissipative and uniformly persistent; in other words, all solutions $x(t;\var),\var\in S$, are defined on $[0,\infty)$, and there are positive constants $m,M$ such that, given any $\var\in S$, there exists $t_0=t_0(\var)$  for which
$$m\le x_i(t;\var)\le M,\q  1\le i\le n,\, t\ge t_0.$$
Here, unless stated otherwise, we take $S=BC_0$ as the set of admissible initial conditions. 

\med

As observed above, from (A4) we get $F(0)\ge 0$. The situation $F_i(0)>0$ for one or more components of the vector $F(0)$ is not excluded from our setting, although  $F(0)=0$ is a quite natural condition:  for many  applications of interest, (3.1) corresponds to a population dynamics model, for which zero should be a steady state. If $F(0)>0$, the persistence of (3.1) follows immediately from the quasimonotone condition (Q): with 
$f(0)=F(0)>0$,  the solution $x^0(t):=x(t;0)$ is  nondecreasing in $t\ge 0$ [30, p.~82] and strictly increasing on an interval $[0,\vare]\, (\vare>0)$. From the order-preserving semiflow, this shows that
  $\liminf_{t\to\infty}x_i(t;\var)\ge \liminf_{t\to\infty}x_i^0(t)\ge x_i^0(\vare)>0, 1\le i\le n,$ for any $\var\in BC_+$.

 In the case of bounded delays, there is an extensive literature using the theory of monotone dynamical system to study the persistence of both autonomous and non-autonomous delayed population models. The situation is much more complex for unbounded delays, since some of the usual methods do not work unless additional conditions are imposed: in spite of a careful choice of  an appropriate phase space, as the one set in Section 2,
 it is a rather difficult task to deal with solutions  with initial conditions $\var=(\var_1,\dots,\var_n)$ in the full set $BC_0$, due to the fact that one may have $\dps \min_{1\le i\le n}\inf_{s\le 0} \var_i (s)=0$.  An alternative way to by-pass this difficulty, is to consider only solutions with initial conditions  whose components  are bounded  below and above  on $(-\infty,0]$ by  positive constants.
The  challenge here is to obtain sufficient conditions for the persistence of (3.1) in $BC_0$ --  a situation not often addressed in the literature, unless there are additional conditions on $F,G$ which allow to relate solutions of (1.1) with solutions of an already known permanent system, or the permanence appears as a by-product   of the global attractivity of a positive equilibrium.   

  The main result about persistence is given below.

 \proclaim{Theorem 3.1}. Suppose that the initial value problems (3.1)-(3.2) have unique solutions  defined on $[0,\infty)$.
 If assumptions (A1)-(A5) are satisfied, system (3.1) is persistent (in $BC_0$).

{\it Proof}. For a positive vector $v=(v_1,\dots,v_n)$ as in (A5),  fix any small $\de >0$  such that 
$$F(v)-(B+\de I) v>0.\eq(3.4)$$ 
Consider $c>0$  as in (A3):  for $\phi,\psi\in BC_+$ with $\phi\le \psi$ and $\phi_i(0)=\psi_i(0)$, 
$$0\le F_i(\psi)-F_i(\phi)\le c\|\psi-\phi\|_g,\q 1\le i\le n.$$
By (A1)-(A2), there exists $\vare>0$ such that
for $\phi,\psi\in BC_+,\phi\le \psi $, with $\phi_i(0)=\psi_i(0)$ and $\|\psi-\phi\|_g< \vare$,
$$  0\le G_i(\phi)-G_i(\psi)< \de  /2\q {\rm for}\ i=1,\dots,n.$$

\smal

{\it Step 1}.
Let $\var\in BC_0$ be given.  Choose $M_0>0$ such that 
$${{2\|\var\|_\infty}\over {g(-M_0)}}< \vare\q {\rm and}\q
{{c|v|}\over {g(-M_0)}}< {\de\over 2}\min_{1\le i\le n}v_i,
\eq(3.5)$$
where  $|v|=|v|_\infty$. Note that $F_i(v)-(v_iG_i(v/m)+\de v_i )\to F_i(v)-(G_i(0)+\de )v_i >0$ as $m\to\infty$.
Fix $m_0\in\N$ such that  
$$F_i(v)-v_iG_i\Big ({1\over m}v\Big)-\de v_i>0, \  i=1,\dots,n,m\ge m_0,
\eq(3.6)$$
and
$$\var(s)> {1\over m_0}v\q {\rm on} \q [-M_0,0].$$

As before, write  (3.1) in the form (3.3).
Assumption (A1) implies that   (3.1) is  cooperative.
Together with (3.1), consider the auxiliary FDE
$$x'(t)=f^\de(x_t),
\eqno (3.7)$$
where $f^\de=(f_1^\de,\dots, f_n^\de), f_i^\de(\phi)=f_i(\phi)-\de \phi_i(0)$ for $\phi\in BC_+,i=1,\dots,n.$
Next, we denote  $x(t)=x(t;\var,f)$ and $y^m(t)=x(t;{1\over m}v,f^\de)$
for any $m\ge m_0$.

 From (3.6) and the sublinearity of $F$, for $m\ge m_0$ and $i=1,\dots,n$,
$$f_i^\de \Big ({1\over m}v\Big )\ge {1\over m}\left [F_i(v)-v_iG_i\Big ({1\over m}v\Big)-\de v_i\right ]>0.$$
Clearly, system (3.7) is cooperative as well, consequently $y_i^m(t)$ is nondecreasing on $[0,\infty)$ ($1\le i\le n$)    (see Theorem 5.1.1 and  Corollary 5.2.2 of [30]).  In particular,   $\inf_{t\ge 0} y_i^m(t)={v_i\over m}>0$.

\med

{\it Step 2}.  We now prove the following claim:
$$x(t)\ge y^m(t)\q {\rm for}\q t\ge 0,m\ge m_0.
\eq(3.8)$$

If the claim fails to be true, there are $m\ge m_0$  and  $t_0>0$, $i\in \{1,\dots,n\}$ such that
$$\eqalign{
&x_j(t)>y_j^m(t)\q {\rm for}\q t\in [0,t_0), j=1,\dots,n,\cr
&x_i(t_0)=y_i^m(t_0).\cr}
$$
Hence, 
$$\eqalign{
&0\ge x_i'(t_0)- (y_i^m)'(t_0)\cr
&=F_i(x_{t_0})-F_i(y_{t_0}^m)+x_i(t_0)[\de-G_i(x_{t_0})+G_i(y_{t_0}^m)]\cr
&=F_i(x_{t_0})-F_i(\tilde \phi)+F_i(\tilde \phi)-F_i(y_{t_0}^m)+x_i(t_0)\left [\de -G_i(x_{t_0})+G_i(\tilde \psi)-G_i(\tilde \psi)+G_i(y^m_{t_0})\right ],\cr}
\eq(3.9)$$
where we take  
$$\tilde\phi(s)=(\tilde\phi_1,\dots,\tilde \phi_n),\q \displaystyle \tilde\phi_i(s)=\system{
&\min (\var_i(t_0+s),{1\over m}v_i),&\ s\le -(M_0+t_0)\cr
&y_i^m(t_0+s),&\  -(M_0+t_0)\le s\le 0\cr
}$$ 
and $$\displaystyle \tilde\psi(s)=\system{
\var (-M_0),&\q  s\le -(M_0+t_0)\cr
x(t_0+s),&\q  -(M_0+t_0)\le s\le 0\cr
}.$$
Since $x_{t_0}(s)\ge \tilde \phi(s),\tilde \psi(s)\ge y^m_{t_0}(s)$ on $(-\infty,0]$ and $\tilde \phi_i(0)=\tilde \psi_i(0)=x_i(t_0)=y^m_i(t_0)$, by (A1) we have $F_i(x_{t_0})-F_i(\tilde \phi)\ge 0,G_i(\tilde \psi)-G_i(y^m_{t_0})\le 0$. On the other hand, from (3.5)
$$\|y^m_{t_0}-\tilde \phi\|_g\le{1\over m}\max_{1\le i\le n}{{v_i}\over{g(-(M_0+t_0))}}\le {{|v|}\over {mg(-M_0)}}<{{\de v_i}\over {2mc}},$$
and
$$\|x_{t_0}-\tilde \psi\|_g=\sup_{s\le -(M_0+t_0)}{{|\var (s+t_0)-\var (-M_0)|}\over{g(s)}}\le {{2\|\var\|_\infty}\over {g(-M_0)}}<\vare.$$
This  implies that $0\le F_i(y_{t_0}^m)-  F_i(\tilde \phi)\le \de v_i/(2m)$ and $|G_i(\tilde \psi)-G_i(x_{t_0})|< \de/2$. From (3.9), we therefore obtain 
$$0>-\de{{v_i}\over{2m}}+{\de\over 2} x_i(t_0).$$
But   $x_i(t_0)=y_i(t_0)\ge{v_i\over m}$ from Step 1. The above inequality yields $0>-\de{{v_i}\over{2m}}+{\de\over 2} x_i(t_0)\ge 0$, which is not possible. This proves the claim (3.8).

\smal
From   Steps 1 and 2,
$$\liminf_{t\to\infty}x_i(t)\ge \liminf_{t\to\infty}y_i^{m_0}(t)\ge {{v_i}\over {m_0}},\q i=1,\dots,n,
\eq(3.10)$$
which shows the persistence of (3.1) in $BC_0$. \ter

\med

{\bf Remark 3.1}. Consider an FDE with {\it bounded} delays, written in abstract form as (3.1) in the usual phase space $C=C([-\tau,0];\R^n)$ with the supremum norm. Let $C_+$ be the  non-negative cone of $C$,  and take $int (C_+)=\{ \var\in C: \var (s)>0$ for $s\in [-\tau,0]\}$ as the set of admissible initial conditions. In this case, we can compare directly the solutions $x(t;\var)$ and $x(t;{1\over m}v)$ of  (3.1), for $m\in\N$ small enough so that $\var(s)\ge {1\over m}v$ on $[-\tau,0]$. In this way,  the persistence of (3.1) in $int (C_+)$ is trivially obtained by assuming (A1), (A4) and (A5) (or, alternatively, (A1) and (A5')). Since the system is autonomous, this leads to the persistence of (3.1) in $C_0:=\{ \var\in C_+: \var (0)>0\}$.

\proclaim{Corollary 3.1}. Consider $F:BC_+\to \R^n$ satisfying (Q), (A3), (A4). If there is a  vector $v>0$ such that $F(v)>0$, the system
$$x'(t)=F(x_t)
\eq(3.11)$$
 is persistent in $BC_0$.

 In the case of finite delays, again the hypothesis (A3) is not required in the above corollary. It should be stressed that, even in the case of finite delays,  this corollary provides a better criterion for persistence of cooperative FDEs  than the one in [35]. 

%
 
 To address the permanence of (3.1), we start by establishing the boundedness of  all solutions.
 
 \proclaim{Lemma 3.2}. Assume (A1)-(A4) and (A6). Then the solutions of the systems (3.1)-(3.2) are defined and bounded on $[0,\infty)$.
 
 {\it Proof}. Again, write (3.1) in the form (3.3). From the sublinearity of $F\big |_{\R^n_+}$, for $q=(q_1,\dots,q_n)>0$ as in (A6) and $L\ge 1$ we have $F(Lq)\le LF(q)$, thus
$$
f_i(Lq)\le L\big (F_i(q)-q_i G_i(Lq)\big)<0,\  i=1,\dots,n.
\eq(3.12)$$
From the quasimonotone condition,  $x_i(t;Lq,f)$ is nonincreasing in $t\in [0,\infty)$
and $x(t;\var,f)\le x(t;Lq,f)$ if $\var\le Lq$ [30]. This proves that all solutions are defined and bounded on $[0,\infty)$.\ter

The main result about permanence  of system (3.1) in the full set $BC_0$ is stated below and requires a fading memory space as a phase space.

 \proclaim{Theorem 3.2}. Let $F,G$ be continuous, bounded on bounded sets of $BC$, and sufficiently regular so that the initial value problems (3.1)-(3.2) have unique solutions. If (A1)-(A6) are satisfied,  system (3.1) is permanent in $BC_0$. Moreover, there are positive  equilibria $x^*,y^*$ of (3.1) such that any solution $x(t)$ of (3.1) with initial condition $x_0=\var\in BC_0$ satisfies
 $$x^*\le \liminf_{t\to\infty} x(t)\le \limsup_{t\to\infty} x(t)\le y^*.
 \eq(3.13)$$

 {\it Proof}. Let $q\in\R^n_+$ be as (A6).
Since bounded positive orbits are precompact in ${\cal B}=C_g^0$ (or ${\cal B}=UC_\ga$ for some $\ga>0$), then $x(t;Lq,f)\searrow y^*(L)$ as $t\to\infty$, where $y^*(L)$ is an equilibrium of (3.1) [30, p.~82].

In a similar way, for $v=(v_1,\dots,v_n)>0$ as in (A5) and $\vare\in (0,1)$ sufficiently small, we get
$$
f_i(\vare v)\ge \vare \big (F_i(v)- v_i G_i(\vare v)\big)>0, \  i=1,\dots,n.
\eq(3.14)
$$
Thus $x(t;\vare v,f)$ is nondecreasing in $t\in [0,\infty)$.
By Lemma 3.2, all solutions with initial conditions in $BC_0$ are bounded,  hence we conclude that  $x(t;\vare v,f)\nearrow x^*(\vare)$ as $t\to\infty$, where $x^*(\vare)$ is a positive  equilibrium of (3.1).

Let $\var\in BC_0$ be given, and take $\de>0$ such that  (3.4) holds.
From Theorem 3.1 (see (3.8)), we also obtain that there exists $m_0\in\N$ with
$x(t;\var,f)\ge x(t;{1\over m}v,f^\de)$ for $t\ge 0$ and $m\ge m_0$.  We further prove the following claims.

\med

{\bf Claim 1}. There is $\vare_0>0$ such that  for $0<\vare\le \vare_0$,
$$x(t;\vare v,f)\nearrow x^*\q {\rm as}\q t\to\infty,
\eq(3.15)$$
where $x^*$ is an equilibrium of (3.1).

\med

To prove this claim, we adapt some arguments in [30, pp.~62]. For similar ideas, see also [29,35].
As shown above, if $\vare_0>0$ is such that $f_i(\vare v)>0$ for  $0<\vare\le \vare_0$ and $1\le i\le n$, then
 $x(t;\vare v,f)\nearrow x^*(\vare)$ as $t\to\infty$, where $x^*(\vare)$ are equilibria, and $x^*(\vare)\le x^*=:x^*(\vare_0)$. Moreover, since $(dx/dt)(0;\vare v,f)=f(\vare v)>0$, then  $x(t;\vare v,f)>x(0;\vare v,f)=\vare v$ for all $t>0$ and $0<\vare\le \vare_0$. Clearly, there is $l^*>0$ such that $x^*\le l^*v$, so  the suprema
$$\de(\vare)=\sup \{ \de \ge 0: x^*(\vare)\ge \de v\}$$
are well defined for  $\vare\in (0,\vare_0]$. 
Note that $\de(\vare)\ge \vare>0$. If $\de(\vare)<\vare_0$, then $x^*(\vare)\ge \de (\vare) v$ implies that
$x^*(\vare)=x(t;x^*(\vare),f)\ge x(t;\de (\vare) v,f)>\de (\vare) v$, which contradicts the definition of $\de(\vare)$. Hence, $\de(\vare)\ge \vare_0$.  We now obtain
$x^*(\vare)=x(t;x^*(\vare),f)\ge x(t;\vare_0 v,f)\nearrow x^*$ as $t\to\infty$. Consequently, $x^*(\vare)=x^*$ for $0<\vare\le \vare_0$.

\med

{\bf Claim 2}. There is $L_0>0$ such that  for $L\ge L_0$,
$$x(t;Lq,f)\searrow y^*\q {\rm as}\q t\to\infty,
\eq(3.16)$$
where $y^*$ is an equilibrium of (3.1).

\med

We argue as in  the proof of Claim 1, so some details are omitted. For $L_0>0$ such that $L_0q>\vare_0v$ and $f_i(Lq)<0$ for $L\ge L_0,1\le i\le n$, we derive $x(t;Lq,f)\searrow y^*(L)$, with 
$Lq\ge y^*(L)\ge y^*$ where $y^*:=y^*(L_0)\ge x^*$. For $L\ge L_0$, define
$$M(L):=\inf \{ M>0:y^*(L)\le Mq\}>0.$$
If $M(L)>L_0$, we get $y^*(L)\le x(t;M(L)q,f)<M(L)q$ for $t>0$, which is not possible from the definition of $M(L)$. Hence $M(L)\le L_0$, and we get $y^*(L)\le x(t;L_0q,f)\to y^*$ as $t\to\infty$. This proves that $y^*(L)= y^*$ for $L\ge L_0$.

\med

{\bf Claim 3}.  For $x^*=x^*(\vare_0)$ and $y^*=y^*(L_0)$ as in Claims 1 and 2, the estimates (3.13) hold; in particular,
(3.1) is permanent in $BC_0$.

\med

To prove the estimates (3.13),  fix any $\var \in BC_0$ and choose $L>L_0$ such that $\var \le Lq$, where $L_0$ is as above. From  Claim 2 it follows that $\limsup_{t\to\infty} x_i(t;\var,f)\le \lim_{t\to\infty} x_i(t;Lv,f)= y_i^*,\ i=1,\dots,n$, hence (3.1) is dissipative.

Next, we remark  that system (3.7), obtained from (3.3) by replacing $G_i(\phi)$ by $G_i(\phi)+\de$ for $1 \le i\le n$, also satisfies the assumptions (A1)-(A6) for any $\de >0$  such that (3.4) holds.
Proceeding as in Claim 1, for $\de >0$ small,
we get that there exists  $\vare_0(\de)>0$ such that  for $0<\vare\le \vare_0(\de)$,
$$x(t;\vare v,f^\de)\nearrow x^*_\de\q {\rm as}\q t\to\infty,$$
where $x^*_\de>0$ is an equilibrium of (3.7). Now, from (3.8), there exists $m_0(\de)\in\N$ such that for 
 $m\ge m_0(\de)$ we have $x(t;\var, f)\ge x(t;{1\over m} v,f^\de)=:y^{m,\de}(t)$. With $m\ge 1/\vare_0(\de)$, we get that 
 $\liminf_{t\to\infty}x(t;\var, f)\ge x^*_\de$.  
 We choose e.g. $\de_k=1/k$ for $k\in\N$, and note that $f^{\de_k}$ increases as $k$ increases. By the definition of $m_0=m_0(\de)$ in the proof of Theorem 3.1 and the definition of $\vare_0=\vare_0(\de)$ in Claim 1, we may take $m_k:=m_0(\de_k)\nearrow\infty$.
 Given $k_1,k_2\in\N$ with $k_2>k_1$, for $m\in\N$ with $m\ge \max (\vare_0(\de_{k_1}),\vare_0(\de_{k_2}))$, we get $y^{m,\de_{k_2}}\ge y^{m,\de_{k_1}}$. Thus we deduce that the sequence of vectors $X^k:=x_{\de_k}^*$ increases at $k$ increases, and is bounded from above by $y^*$, thus there exists $x^*:=\lim_k X_k^*$. Moreover, $X^k=(X^k_1,\dots,X^k_n)$ is an equilibrium of (3.7) with $\de=1/k$, i.e., $f_i(X^k)-{1\over k} X^k_i=0$, and $\liminf_{t\to\infty}x(t;\var, f)\ge X^k$. By letting $k\to\infty$, we obtain the estimate
$\liminf_{t\to\infty}x(t;\var, f)\ge x^*$, where $x^*$ is an equilibrium of (3.1).\ter

\proclaim{Corollary 3.2}. Under conditions (A1)-(A6), there exists a positive equilibrium of (3.1). Furthermore, if this positive equilibrium is unique, then it is a global attractor of all solutions with initial conditions in $BC_0$.

 {\bf Remark 3.2}. As a consequence of the above  proof,  the set $\{\var\in BC: \vare_0v\le \var\le L_0q\}$ is positively invariant and a global attractor for the semiflow of (3.1). However, if the equilibria $x^*,y^*$ in (3.13) are distinct, this set may contain a very complex dynamics, with more equilibria, periodic orbits, or heteroclinic orbits (see e.g. [22,23]).

\smal

The next result is important for applications.

\proclaim{Corollary 3.3}. Consider an FDE with infinite delay of the form
$$
x_i'(t)=F_i(x_t)-x_i(t)G_i(x_i(t)), \quad i=1,\dots,n
\eq(3.17)$$
where    $F_i:BC\to\R, G_i:[0,\infty)\to \R$ are   Lipschitz continuous on bounded  sets of their domains,  and $G_i(0)=0$, $1\le i\le n$. Suppose that $F$ satisfies (Q), (A3), (A4), and that there are positive vectors $v, q$ such that $F(v)>0,F(q)<0$. Then the solutions of (3.17) with initial conditions $x_0=\var\in BC_0$ are defined on $[0,\infty)$ and (3.17) is permanent in $BC_0$.

{\it Proof}. Since $F_i,G_i$ are Lipschitz continuous on bounded sets, initial value problems have unique non-conti\-nuable solutions.
For (3.17), conditions (A1)(ii) and (A2) are trivially satisfied,  and  $F(q)<0$ for some vector $q>0$ implies (A6).  The latter also implies  that solutions are defined on $[0,\infty)$. The permanence follows from Theorem 3.2.\ter

\bigskip

{\bf 4. Persistence and permanence for non-autonomous systems}

\med
We now consider a class of non-autonomous FDEs with infinite delay given by
$$
x_i'(t)=F_i(t, x_t)-x_i(t)G_i(t,x_t),\quad i=1,\dots,n,
\eq(4.1)$$
for $t\ge t_0$, where    $F,G:D\subset \R \times C_g^0\to\R^n$ are  continuous and bounded on bounded sets of $D$,
$a\in\R, D=(a,\infty)\times \Omega$, $t_0>a$ and $ \Omega$  is an open set containing $BC_+$.
Initial conditions are  written as
$x_{t_0}=\var$ with $\var\in BC_0.$  In this section, we always assume that solutions to FDEs with initial conditions in $BC_0$ are unique and defined on $[0,\infty)$.

\proclaim{Theorem 4.1}. (i) Assume that
there are continuous functions $F^l, G^u:BC\to \R^n$  such that
$$
F^l(\phi)\le F(t,\phi),\ G(t,\phi)\le G^u (\phi)\q {\rm for}\q (t,\phi)\in D.
\eq(4.2)$$
If the pair $(F^l,G^u)$ satisfies (A1)-(A5),  (4.1) is persistent in $BC_0$.\vskip 0cm
(ii) Assume that
there are continuous functions $F^l,F^u, G^l, G^u:BC\to \R^n$  such that
$$
F^l(\phi)\le F(t,\phi)\le F^u (\phi), \ G^l(\phi)\le G(t,\phi)\le G^u (\phi)\q {\rm for}\q (t,\phi)\in D.
\eq(4.3)$$
 If  $(F^l,G^u)$ satisfies (A1)-(A5) and $(F^u, G^l)$ satisfies (A1)-(A6),  (4.1) is permanent in $BC_0$. 

{Proof}. 
In what follows, and without loss of generality,  suppose that $t_0=0$.
Write (4.1) as 
$$x'(t)=f(t,x_t),\ t\ge 0,\eq(4.4)$$ where $f_i(t,\phi)=F_i(t,\phi)-\phi_i(0)G_i(t,\phi), 1\le i\le n$.  
 Under (4.2), consider the autonomous system
$$
x_i'(t)=F_i^l(x_t)-x_i(t)G_i^u(x_t)=:f_i^l(x_t), \quad i=1,\dots,n,
\eq(4.5)$$
and assume that $(F^l,G^u)$ satisfies (A1)-(A5). Since (Q) holds for $f^l$, for non-negative solutions of (4.1) we have  $x(t;0,\var, f)\ge x(t;\var,f^l)$ ($\var\in BC_0$). Theorem 3.1 guarantees that (4.5) is persistent, so (4.1) is persistent as well.
Similarly, if (4.3) holds we further compare the solutions of (4.1) with the solutions of the cooperative system
$$
x_i'(t)=F_i^u(x_t)-x_i(t)G_i^l(x_t)=:f^u_i(x_t), \quad i=1,\dots,n,
\eq(4.6)$$
Since  $x(t;\var,f^l)\le x(t;0,\var, f)\le x(t;\var,f^u)$  for $\var\in BC_0$,  (4.1) is permanent.\ter

Under conditions (4.3),  in fact we obtain  explicit lower and upper uniform bounds for  the solutions $x(t;0,\var, f),\var\in BC_0$, of (4.1):
 $$x^{*,l}\le \liminf_{t\to\infty} x(t;0,\var, f)\le \limsup_{t\to\infty} x(t;0,\var, f)\le y^{*,u},
$$
 where $x^{*,l},y^{*,u}$ are equilibria of (4.5), (4.6), respectively.

 Theorem 4.1 is not easily applicable to FDEs with unbounded time-dependent discrete delays, in which case it  is better to deal directly with the non-autonomous system.

 \proclaim{Theorem 4.2}. For (4.1),  let the pair $(F,G)$ satisfy the assumptions:\vskip 1mm
 {\parindent=1.5cm{
\item{(H1)} 
 for $(t,\phi),(t,\psi)\in D$ with $0\le \phi\le \psi,\phi_i(0)=\psi_i(0)$, then (i) $F_i(t,\phi)\le F_i(t,\psi)$ and (ii) $G_i(t,\phi)\ge G_i(t,\psi), \ i=1,\dots,n$;
   \item{(H2)} for each $\de>0$ and each $ i=1,\dots,n$, there exists  $\vare>0$ such that for all $t\ge t_0, \phi,\psi\in BC_+,\phi\le \psi$ with $\|\psi-\phi\|_g<\vare$ and $\phi_i(0)=\psi_i(0)$, then $G_i(t,\phi)-G_i(t,\psi)<\de$;
\item{(H3)} there exists a constant $c>0$ such that for $t\ge t_0,\phi,\psi\in BC_+,\phi\le \psi$ and $\phi_i(0)=\psi_i(0)$, then $F_i(t,\phi)-F_i(t,\psi)\ge -c\|\psi-\phi\|_g, \ i=1,\dots,n$.
}}
\vskip 1mm
In addition, suppose that  the functions $x\mapsto F(t,x)=:\hat F(x)$ and $x\mapsto G(t,x)=:\hat G(x), x\in\R^n,$ do not depend on $t$. Then,\vskip 0cm
(i) if the pair $(\hat F, \hat G)$ satisfies (A4) and (A5), system (4.1) is persistent in $BC_0$.\vskip 0cm
(ii) if the pair $(\hat F, \hat G)$ satisfies (A4) and (A6), all solutions of (4.1) with initial conditions in $BC_0$ are bounded.

{\it Proof}.  For non-autonomous systems  $x'(t)=f(t,x_t)$
 with $f:D\subset \R\times C_g^0\to\R^n$ continuous, the quasimonotone  condition (Q) should be replaced by

\vskip 3mm

\item{(Q')} for $(t,\phi), (t,\psi)\in D$ with $0\le \phi\le \psi$ and $\phi_i(0)=\psi_i(0)$, then $f_i(t,\phi)\le f_i(t,\psi), \ i=1,\dots,n$.

\vskip 3mm
{\parindent=0cm Condition (Q') implies that the solution operator $x_t(\sigma,\phi)$ for $x'(t)=f(t,x_t)$ is monotone relative to the variable $\phi$ and that comparison results as in  Chapter 5 of [30] are valid. }

(i) As before, write (4.1) as (4.4) and take e.g. $t_0=0$.
Assuming (H1)-(H3) means that (A1)-(A3) are satisfied with $F(\cdot)$ and $G(\cdot)$ replaced by $F(t,\cdot)$ and $G(t,\cdot)$, respectively ($t\ge 0)$. In particular,  (H1) implies  (Q'). 

We now adapt  the proof of 
Theorem 3.1, by replacing equation (3.1)  by (4.1) and (3.7) by
$$x'(t)=f^\de (t,x_t),\eq(4.7)$$
where $f_i^\de(t,\phi)=f_i(t,\phi)-\de \phi_i(0),1\le i\le n$, and $\de>0$ is such that $\hat F(v)-(B+\de I)v>0$ for $v=(v_1,\dots,v_n)>0$ and $B=diag\, (\hat G_1(0),\dots,\hat G_n(0))$ as in (A5).  Clearly, (4.7) also satisfies (Q').

Let $\var\in BC_0$ be given. It is apparent that the arguments in Step 2 of the proof of 
Theorem 3.1  apply  with very little change to the non-autonomous case, thus there exists $m_0\in\N$ such that $x(t)\ge y^m(t)$ for $t\ge 0,m\ge m_0$, where $x(t):=x(t;0,\var, f)$ and  $y^m(t):=x(t;0,{1\over m},f^\de)$ . 

To further prove that $\inf_{t\ge 0}y_i^m(t)>0$ for all $i$ and $m$ large, we replace Step 1 of the proof of Theorem 3.1 by the following argument: for $t\ge 0$ and $\phi\ge {1\over m}v$ with  $\phi_i(0)={1\over m}v_i$, from (H1) and (A4) we have 
$$f_i^\de(t,\phi)\ge F_i(t,v/m)-{1\over m}v_i\Big (G_i(t,v/m)+\de\Big)\ge {1\over m}\Big [\hat F_i(v)-v_i\big (\hat G_i(v/m)+\de\big )\Big]>0$$
for $m\in \N$ sufficiently large and $1\le i\le n$; from Remark 5.2.1  of Smith's monograph [30], we derive that $y^m(t)\ge  {1\over m}v$ for $t\ge 0$. This ends the proof of the assertion (i).

(ii)  In a similar way,  under the hypotheses in (ii), a few changes to the proof of  Lemma 3.2 allow us  to conclude that all solutions are bounded.\ter
%

The permanence result in Theorem 3.2 is not easily adapted to deal directly with non-autono\-mous systems, since its proof uses  $\omega$-limit sets for autonomous FDEs. In fact, an argument often used in the proof is that a bounded solution $x(t)$ with monotone components should converge to some $x^*$ as $t\to\infty$, where $x^*$ is an equilibrium of the system -- which obviously  need not  happen for a non-autonomous system. 
 However, in concrete applications, one might further derive the permanence of the system, with explicit upper and lower bounds, as shown in the next  section.

\bigskip

{\bf 5. Applications to  population models}

\med

The above results on persistence and permanence apply to a very broad class of abstract FDEs with infinite delay, which include many important models used in population dynamics. An advantage of our method is that, in general,  the validity of hypotheses (A1)-(A6) and (H1)-(H3)   is  very easy to check. Here, we illustrate the results with some selected examples.


\med

{\bf 5.1. A cooperative Lotka-Volterra model with patch structure}
\med

 Consider the following cooperative Lotka-Volterra system with patch structure:
$$\eqalign{
x_i'(t)=x_i(t)\biggl(\be_i-\mu_i x_i(t)&+\sum_{j=1}^n a_{ij} \int_0^{\infty}x_j(t-s)\, d\eta_{ij}(s)\biggr)\cr
&+\sum_{j=1}^n d_{ij}\int_0^{\infty}x_j(t-s)\, d\nu_{ij}(s), \quad i=1,\ldots,n\cr}
\eq(5.1)$$
where   $\be_i, \mu_i\in\R,a_{ij}\ge 0, d_{ij}\ge 0$,  $\eta_{ij},\nu_{ij}:[0,\infty)\to \R$  are bounded, nondecreasing functions with total variation one, for all $i,j$.  For results on persistence, permanence and stability for autonomous or non-autonomous Lotka-Volterra with infinite delays, as well as for a biological explanation of the coefficients involved, see [8,10,11,20,21,24,32,33]  and references therein.
System (5.1) is written in the form (3.1) with $F_i, G_i:BC\to\R$ linear  operators given by
$$\eqalign{
F_i(\phi)&=\be_i\phi_i(0)+\sum_{j=1}^n d_{ij}\int_0^{\infty}\phi_j(-s)\, d\nu_{ij}(s),\cr
G_i(\phi)&=\mu_i \phi_i(0)-\sum_{j=1}^n a_{ij} \int_0^{\infty}\phi_j(-s)\, d\eta_{ij}(s).\cr
}\eq(5.2)$$
We now insert (5.1) in a phase space $C_g^0$,  to justify the use of Theorems 3.1 and 3.2.
As shown in [12,15], given any prescribed $\de >0$, there is a function $g$ satisfying (g1)-(g3) and such that
$$\int_0^{\infty} g(-s)\, d\nu_{ij}(s)<1+\de,\q \int_0^{\infty}g(-s)\, d\eta_{ij}(s)<1+\de,\q i,j=1,\dots,n.$$ 
In this way, $F_i,G_i$ are well-defined by (5.2) as bounded linear operators in $C_g^0$ [8,18], whose operator norms satisfy the estimates
$\|F_i\|<|\be_i|+(1+\de)\sum_{j=1}^n d_{ij},\, \|G_i\|<|\mu_i|+(1+\de)\sum_{j=1}^n a_{ij},\, i,j=1,\dots,n.$
Clearly, assumptions (A1)-(A4) are fulfilled. Now define the matrices
$$M=diag\, (\be_1,\dots,\be_n)+[d_{ij}],\q N=diag\, (\mu_1,\dots,\mu_n)-[d_{ij}].$$
For any $v\in\R^n$, $F(v)=Mv, G(v)=Nv$. For the present setting, we derive the following criteria:

\proclaim{Theorem 5.1}. With the above notations,\vskip 0cm
 (i) If  there exists a positive vector $v$ such that $Mv>0$,  system (5.1) is persistent in $BC_0$. \vskip 0cm (ii) If  there exist positive vectors $v$ and $q$ such that $Mv>0$ and $Nq>0$,  system (5.1) is permanent in $BC_0$ and it has a positive equilibrium.\vskip 0cm
 (iii) In  the latter case, if there is a positive equilibrium $x^*$ for which $Mx^*>0$, then $x^*$ is a global attractor for all  solutions of (5.1) with initial conditions in $BC_0$.
  
  {\it Proof}. If $Mv>0$ for some positive $v\in\R^n$, then (A5) is satisfied; similarly,  $Nq>0$ for some  $q>0$ implies  (A6). The assertions  (i) and (ii) follow from Theorems 3.1 and 3.2.
  
  We now further assume that $Mx^*>0$ for some positive equilibrium. In particular, this means that $Mv>0$ for $v=\vare x^*$ and $\vare >0$;  moreover, for $\vare$ small, say $\vare\in (0,1)$,    as in (3.15) we have  $x(t;\vare x^*)\nearrow x^*$ as $t\to\infty$. Then  any solution $x(t)=x(t,\var)\ (\var\in BC_0^+)$ of (5.1) satisfies the estimates (3.13), where $y^*$ is also an equilibrium.  We need to show  that $x^*=y^*$. 
  
   In order to simplify the arguments, we consider the ODE system associated with (5.1), obtained by taking all the delays equal to zero:
 $$
x_i'(t)=x_i(t)\biggl(\be_i-\mu_i x_i(t)+\sum_{j=1}^n a_{ij} x_j(t)\biggr)+\sum_{j=1}^n d_{ij}x_j(t), \q i=1,\dots,n.
\eq(5.3)$$ 
Since the equilibria of (5.1) and (5.3) are the same, we only need to show that solutions of (5.3) satisfy  $\limsup_{t\to\infty} x(t)\le x^*$. For a positive solution $x(t)$ of (5.3), define
$L_j=\limsup_{t\to\infty} (x_j(t)/x_j^*)$, and choose $i$ such that $L_i=\max_{j}L_j$. Choose a sequence $t_k\to \infty$ with $x_i(t_k)\to L_ix_i^*$ and $x_i'(t_k)\to 0$ as $k\to\infty$. Using (5.3), this yields
$$\eqalign{
0&\le L_ix_i^*\biggl(\be_i-\mu_iL_ix_i^*+L_i\sum_{j=1}^n a_{ij} x_j^*\biggr)+L_i\sum_{j=1}^n d_{ij}x_j^*\cr
&=L_i(\be_ix_i^*+\sum_{j=1}^n d_{ij}x_j^*)+L_i^2x_i^*(-\mu_ix_i^*+\sum_{j=1}^n a_{ij} x_j^*)\cr 
&=L_i(1-L_i)(\be_ix_i^*+\sum_{j=1}^n d_{ij}x_j^*)=L_i(1-L_i)(Mx^*)_i.\cr}$$
Under the condition $Mx^*>0$, this leads to  $L_i\le 1$. From  (3.13) and the definition of $L_i$, we get
$x_j^*\le \liminf_{t\to\infty} x_j(t)\le \limsup_{t\to\infty} x_j(t)\le x_j^*$ for $1\le j\le n.$
  \ter
  
\smal

{\bf Remark 5.1}. As expected, Theorem 5.1 shows that the global dynamics of (5.1) depend heavily on the algebraic properties of the matrices $M,N$. Condition $Nq>0$ for some vector $q>0$ is equivalent to saying that $N$ is a non-singular M-matrix; the dissipativeness  follows from this property. When $M$ is an  irreducible matrix,  there is a positive vector $v$ such that $Mv>0$ if and only if $s(M)>0$, where $s(M)$ is the {\it spectral bound} of $M$, defined by $s(M)=\max \{ {\rm Re}\, \la: \la\in\sigma (M)\}$. If all the coefficients $\be_i$ are positive, clearly   $Mx^*>0$ if $x^*>0$ -- so, if $Nq>0$ for some $q>0$, a positive equilibrium  of (5.1) must be a global attractor.
Theorem 5.1  significantly improves  the results in [10] in what concerns the persistence of cooperative Lotka-Volterra models (5.1): here, the persistence is established in $BC_0$, a question which was left as an open problem in [10].
On the other hand,  [10] deals with questions of extinction, persistence, global asymptotic stability for Lotka-Volterra systems (5.1) without any prescribed 
signs for the coefficients $\be_i,a_{ij}$. If there is at least one negative   coefficient $a_{ij}$, then (5.1) is not cooperative, and the theory of monotone systems cannot be applied directly. In this scenario, auxiliary cooperative FDEs can be used, and results of comparison with cooperative systems  invoked, to draw conclusions about the non-cooperative system -- a technique often exploited in [10]. This technique    will be illustrated in Subsection 5.3. 

 \smal
 
 Consider now the following non-autonomous version of the Lotka-Volterra system (5.1):$$\eqalign{
x_i'(t)=x_i(t)\biggl(\be_i(t)-\mu_i(t) x_i(t)&+\sum_{j=1}^n a_{ij} (t)\int_0^{\infty}x_j(t-s)\, d\eta_{ij}(s)\biggr)\cr
&+\sum_{j=1}^n d_{ij}(t)\int_0^{\infty}x_j(t-s)\, d\nu_{ij}(s), \quad t\ge 0,i=1,\ldots,n\cr}
\eq(5.4)$$
where   $\mu_i(t),\be_i(t),a_{ij}(t), d_{ij}(t)$ are continuous and bounded on $[0,\infty)$,  $a_{ij}(t), d_{ij}(t)$ are non-negative and $\eta_{ij},\nu_{ij}:[0,\infty)\to \R$  are bounded, nondecreasing functions with total variation one, for all $i,j$. A straightforward application of Theorems 4.1 and 5.1 gives the following criterion:
 
 \proclaim{Theorem 5.2}. Under the above conditions, define the $n\times n$ constant matrices
 $$\displaylines{M^l=diag\, (\ul \be_1,\dots, \ul \be_n)+\Big[\ul d_{ij}\Big ], \ M^u=diag\, (\ol \be_1,\dots, \ol \be_n)+\Big[\ol d_{ij}\Big ],\cr
N^l=diag\, (\ul \mu_1,\dots, \ul \mu_n)-\Big[\ol a_{ij}\Big ],}$$
where we use the notations
$$\underline f=\inf_{t\ge 0}f(t),\q \ol f=\sup_{t\ge 0}f(t),
\eq(5.5)$$
for a function $f:[0,\infty)\to\R $ bounded. Then:
 (i) if  there exists a  vector $v>0$ such that $M^lv>0$,  system (5.4) is persistent in $BC_0$; (ii) if  there exist  vectors $v>0$ and $q>0$ such that $M ^uv>0$ and $N^lq>0$,  system (5.4) is permanent in $BC_0$.
 
\

  {\bf 5.2. A system of modified logistic equations with unbounded delays}
  
 \med
 Consider the classical delayed logistic equation, also known as Wright's equation, given by $x'(t)=rx(t)(1-x(t-\tau)/K)$, where $x(t)$ is the population of a species at time $t$, $r$ the intrinsic growth rate of the species, $K$ the carrying capacity, and $\tau$ the maturation delay.
  In Arino et al. [6], the model
 $$x'(t)={{\gamma \mu x(t-\tau)}\over {\mu e^{\mu \tau}+k (e^{\mu \tau}-1)x(t-\tau)}}-\mu x(t)-\k x^2(t),
\eq(5.6)$$
where $\gamma,\mu,\k ,\tau >0$, was derived as an alternative, and more realistic, formulation for the delayed logistic  law.
 The coefficients in the  logistic equation are related to the ones in (5.6)   by $r=\gamma -\mu$  (for $\gamma, \mu$ the birth and mortality rates, respectively) and $K=(\gamma -\mu)/\k$.
More recently, Bastinec et al. [7] proposed a more general {\it non-autonomous} model with multiple time dependent delays:
 $$
 x'(t)=\sum_{k=1}^m {{\al_k(t)x(t-\tau_k(t))}\over {1+\be_k(t)x(t-\tau_k(t))}}
-\mu(t)x(t)-\k(t)x^2(t),\q t\ge 0,
\eq(5.7)$$
where  $\al_k,\k:[0,\infty)\to (0,\infty), \be_k,\mu,\tau_k:[0,\infty)\to [0,\tau]$ (for some $\tau>0$) are continuous, $1\le k\le m$. The question of the permanence of (5.7)  was raised (but not studied) in [7],  and a criterion for it given  by the author in [9], with explicit uniform lower and upper bounds for all positive solutions. 
 
 Here, we pursue the research, and generalize the scalar equation (5.7): not only {\it unbounded} time-varying delays are allowed, but also we consider $n$ classes (of one or multiple species) following the `modified delayed logistic equation'  (5.7), with dispersal terms among the classes. This leads to
  the system
$$
\eqalign{
x_i'(t)=\sum_{k=1}^{m_i} {{\al_{ik}(t)x_i(t-\tau_{ik}(t))}\over {1+\be_{ik}(t)x_i(t-\tau_{ik}(t))}}&+\sum_{j=1}^nd_{ij}(t) x_j(t-\sigma_{ij}(t))\cr
&-\mu_i(t)x_i(t)-\k_i(t)x_i^2(t), \quad t\ge 0,\, i=1,\dots,n\cr}
\eq(5.8)$$
where all the coefficients $\al_{ik}(t),\be_{ik}(t), d_{ij}(t), \mu_i(t), \k_i(t)$ and delays $\tau_{ik}(t),\sigma_{ij}(t)$  are non-negative and continuous functions in $t\in[0,\infty)$, $k=1,\dots,m_i, i,j=1,\dots,n$; the  time-dependent delay functions $\tau_{ik}(t),\sigma_{ij}(t)$  are  possibly unbounded.  
The migration rates of populations moving from class $j$ to class $i$ are given by $d_{ij}(t)$, with $\sigma_{ij}(t)$  the time-delays during dispersion; the instantaneous loss term $-d_{ji}(t) x_j(t)$ may be incorporated in the term $-\mu_j(t)x_j(t)$ of the $j$th-equation.
For biological reasons, we usually consider $d_{ii}(t)\equiv 0$; in any case for each $i\in \{1,\dots,n\}$ the term $d_{ii}(t) x_i(t-\sigma_{ii}(t))$ may be inclued in the first sum on the right-hand side of (5.8).

We now established sufficient conditions for the persistence of (5.8) in $BC_0$. This example also illustrates how to combine the techniques in Theorems 4.1 and 4.2.

\proclaim{Theorem 5.3}. Suppose that $\al_{ik}, \be_{ik}, d_{ij}, \mu_i, \k_i,\tau_{ik},\sigma_{ij}:[0,\infty)\to [0,\infty)$ are continuous functions,   $A_i(t):=\sum_{k=1}^{m_i}\al_{ik}(t)$ and $ \k_i(t)$ are  bounded  below and above by positive constants, and $ \be_{ik}, d_{ij}, \mu_i$ are  bounded, $k=1,\dots,m_i, i,j=1,\dots,n$. If 
there exists a positive vector $v=(v_1,\dots,v_n)$ such that 
$$Hv>0,\eq(5.9)
$$
 where $H$ is the $n\times n$ matrix
$H= diag\, (\underline A_1-\ol \mu_1, \dots,\underline A_n-\ol \mu_n)+\big [ \underline d_{ij}\big ]$
for $$\underline d_{ij}=\inf_{t\ge 0}d_{ij}(t),\q\underline A_i=\inf_{t\ge 0} A_i(t),\q \ol \mu_i= \sup_{t\ge 0} \mu_i(t),$$
then all solution of (5.8) with initial conditions in $BC_0$ are bounded below and above by positive constants.

{\it Proof}. As in (5.5), we set the notations
$\underline f=\inf_{t\ge 0}f(t),\ \ol f=\sup_{t\ge 0}f(t)$
for a function $f:[0,\infty)\to[0,\infty) $ bounded. 
Condition (5.9) reads as
$$(\underline A_i-\ol \mu_i)v_i+\sum_{j=1}^n \ul d_{ij}v_j>0,\q i=1,\dots,n.
\eq(5.10)$$
Effecting a scaling of the variables, $\hat x_i(t)=x_i(t)/v_i$,
we obtain a new system with the form (5.8), where $\be_{ik}(t), d_{ij}(t),\k_i(t)$ are replaced by $\hat \be_{ik}(t)=v_i\be_{ik}(t), \hat d_{ij}(t)=v_jd_{ij}(t)/v_i,\hat \k_i(t)=v_i\k_i(t)$ for all $i,k,j$, and all  other coefficiens are the same. Dropping the hats for the sake of simplification, we can therefore assume that (5.8) satisfies condition (5.10) with $v=(1,\dots,1)$.

Define $r_{ik}(t,x):={{\al_{ik}(t)x}\over {1+\be_{ik}(t)x}}, t\ge 0, x\ge 0$, and note that $0\le  {{\p}\over {\p x}}r_{ik}(t,x)\le \ol \al_{ik}$. Write (5.8) as 
$$x_i'(t)=F_i(t, x_t)-x_i(t)G_i(t,x_i(t))=:f_i(t,x_t),\quad t\ge 0, i=1,\dots,n,$$ with
$$\eqalign{
F_i(t,\phi)&=\sum_{k=1}^{m_i}r_{ik}(t,\phi_i(-\tau_{ik}(t))+ \sum_{j=1}^nd_{ij}(t) \phi_j(-\sigma_{ij}(t)),\q t\ge 0,\phi\in BC,\cr
G_i(t,x)&=\mu_i(t)+\k_i(t)x, \q t\ge 0, x\in\R.\cr}$$
Define also the functions
$F^l(t,\phi),F^u(t,\phi)$ and $G^l(x),G^u(x)$ whose components are given by
$$\eqalign{
F_i^l(t,\phi)&=\sum_{k=1}^{m_i}{{ \ul\al_{ik}\phi_i(-\tau_{ik}(t))}\over {1+\ol \be_{ik}\phi_i(-\tau_{ik}(t))}}+ \sum_{j=1}^n\ul d_{ij} \phi_j(-\sigma_{ij}(t)),\cr
F_i^u(t,\phi)&=\sum_{k=1}^{m_i}{{ \ol\al_{ik}\phi_i(-\tau_{ik}(t))}\over {1+\ul \be_{ik}\phi_i(-\tau_{ik}(t))}}+ \sum_{j=1}^n\ol d_{ij} \phi_j(-\sigma_{ij}(t)), \q t\ge 0, \phi\in BC,\cr
G_i^l(x)&=\ul \mu_i+\ul \k_ix, \q
G_i^u(x)=\ol \mu_i+\ol \k_ix, \q x\in\R.\cr}$$
Clearly,
$$F_i^l(t,\phi)\le F_i(t,\phi)\le F_i^u(t,\phi),\q G_i^l(x)\le G_i(t,x)\le G_i^u(x)$$ for $ t\ge 0, \phi\in BC_+, x\in\R_+.$
On the other hand, choose $g$ satisfying (g1)-(g3), and insert (5.8) is a space $C_g^0$. The functions
$F^l(t,\phi),F^u(t,\phi)$ are uniformly (for $t\in [0,\infty)$) Lipschitz continuous and nondecreasing with respect to $\phi$ in the non-negative cone $BC_+$, and $G_i^l(x),G_i^u(x)$ are increasing for $x\ge 0$. Therefore  
the pairs $(F^l,G^u)$ and $(F^u,G^l)$ satisfy (H1)-(H3). Note also that $x\mapsto F^l(t,x),x\mapsto F^u(t,x)$  are autonomous for $x\in \R^n$, and that $\hat F^l(x):=F^l(t,x),\hat F^u(x):=F^u(t,x)$ are sublinear.

Next, consider the auxiliary systems
$$\eqalignno{
x_i'(t)&=F_i^l(t,x_t)-x_i(t)G_i^u(x_i(t))=:f^l_i(t,x_t), \quad i=1,\dots,n,
&(5.11)\cr
x_i'(t)&=F_i^u(t,x_t)-x_i(t)G_i^l(x_i(t))=:f^u_i(t,x_t), \quad i=1,\dots,n.
&(5.12)\cr}$$ 
By results of comparison of solutions,
$$x(t;0,\var,f^l)\le x(t;0,\var,f)\le x(t;0,\var,f^u),
\eq(5.13)$$
 for $\var\in BC_0$. From (5.10)   with $v=(1,\dots,1)=:{\bf 1}\in \R^n_+$ and $\vare>0$ small, we have
$$\hat F_i^l(\vare {\bf 1})-G_i^u(0)=\vare \left ( \sum_{k=1}^{m_i}{{ \ul\al_{ik}}\over {1+\vare \ol \be_{ik}}}-\ol \mu_i+\sum_{j=1}^n \ul d_{ij}\right)>0,\q i=1,\dots,n.$$
 Theorem 4.2(i) yields the persistence of (5.11), and thus the persistence of (5.8) as well. Moreover, $\lim_{x\to\infty} G_i^l(x)=\infty,\, 1\le i\le n$, which implies that $(\hat F^u,G^l)$ satisfies (A6). Theorem 4.2(ii) and (5.13) allow us to conclude the boundedness of solutions  $x(t;0,\var,f)$ on $[0,\infty)$, for any $\var\in BC_0$. \ter
 
 With an additional condition on the (possibly unbounded) delay functions, the permanence of (5.8) is established, with explicit lower and upper bounds for the asymptotic behavior of solutions.
 
 \proclaim{Theorem 5.4}. Assume the hypotheses of Theorem 5.3, and in addition that $\lim\limits_{t\to\infty} (t-\tau_{ik}(t))=\lim\limits_{t\to\infty} (t-\sigma_{ij}(t))=\infty$ for $i,j=1,\dots,n, k=1,\dots,m_i$. Then,
 system (5.8) is permanent. Moreover, for $v=(v_1,\dots,v_n)$ as in (5.9), the solutions $x(t)=x(t;0,\var)$ of (5.8) with $\var\in BC_0$ satisfy  the uniform estimates 
$$
m_0\le \liminf_{t\to\infty} (x_i(t)/v_i)\le \limsup_{t\to\infty} (x_i(t)/v_i)\le M_0,\q i=1,\dots,n,
\eq(5.14)$$
where 
$$
M_0=\max_{1\le i\le n}\limsup_{t\to\infty}{1\over {v_i^2\k_i(t)}} \left [v_i\left (\sum_{k=1}^{m_i}\al_{ik}(t)-\mu_i(t)\right )+\sum_{j=1}^n d_{ij}(t)v_j\right]
\eq(5.15)$$
 and
 $$
m_0=\min_{1\le i\le n}\liminf_{t\to\infty}{1\over {v_i^2\k_i(t)}} \left [v_i\left (\sum_{k=1}^{m_i}{{\al_{ik}(t)}\over {1+\be_{ik}(t)v_iM_0}}-\mu_i(t)\right )+\sum_{j=1}^n d_{ij}(t)v_j\right].
\eq(5.16)$$

{\it Proof}. 
 For  $x(t):=x(t;0,\var,f)$, set $ \ul x_j:=\liminf_{t\to\infty} x_j(t), \ol x_j:=\limsup_{t\to\infty} x_j(t)$, $1\le j\le n$. From Theorem 5.3, $0<\ul x_j\le \ol x_j<\infty$ for all $j$.   Consider $i$ such that $\ol x_i/v_i=\max_{1\le j\le n}(\ol x_j/v_j)$, for $v>0$ as in (5.9). By the fluctuation lemma, take a sequence $(t_m)$ with $t_m\to\infty$,
$x_i'(t_m)\to 0$ and $x_i(t_m)\to \ol x_i$. For any $\vare>0$ small and $m$ sufficiently large, we have $x_i(t_m-\tau_{ik}(t_m))\le \ol x_i+\vare$ and  $v_ix_j(t_m-\sigma_{ij}(t_m))/v_j\le \ol x_i+\vare$ for all $k,j$. Recalling that the functions $r_{ik}(t,\cdot)$ are nondecreasing on $x\in [0,\infty)$,  for sufficiently   large $m$ we derive
$$
\eqalign{
x_i'(t_m)\le &(\ol x_i+\vare)\sum_{k=1}^{m_i} {{\al_{ik}(t_m)}\over {1+\be_{ik}(t_m)(\ol x_i+\vare)}}+{{(\ol x_i+\vare)}\over v_i}\sum_{j=1}^nd_{ij}(t_m) v_j -\mu_i(t_m)x_i(t_m)-\k_i(t_m)x_i^2(t_m),\cr
\le & \k_i(t_m)\left [{{\ol x_i+\vare}\over {v_i\k_i(t_m)}}\left (\sum_{k=1}^{m_i} v_i\al_{ik}(t_m)+\sum_{j=1}^nd_{ij}(t_m) v_j\right )-{{\mu_i (t_m)}\over {\k_i(t_m)}}x_i(t_m)-x_i^2(t_m)\right ].\cr
}
$$
Taking limits $m\to\infty,\vare \to 0^+$, this estimate yields
$$0\le \limsup_{t\to\infty}\left [{{\ol x_i}\over {v_i\k_i(t)}}\left (v_i\Big (\sum_{k=1}^{m_i} \al_{ik}(t)-\mu_i(t)\Big )+\sum_{j=1}^nd_{ij}(t) v_j\right )-\ol x_i^2\right ],$$
thus
$$\ol x_i\le \limsup_{t\to\infty}{1\over {v_i\k_i(t)}} \left [v_i\left (\sum_{k=1}^{m_i}\al_{ik}(t)-\mu_i(t)\right )+\sum_{j=1}^n d_{ij}(t)v_j\right].$$
This leads to $\ol x_j/v_j\le M_0$, for $1\le j\le n$ and $M_0$ as in (5.15).

To prove the uniform lower bound given by (5.14), (5.16),  we  reason along the lines above, and some details are omitted. Let $i\in\{1,\dots,n\}$ be such that $\ul x_i/v_i=\min_{1\le j\le n}\ul x_j/v_j$, and
 take a sequence $(s_m)$ with $s_m\to\infty$,
$x_i'(s_m)\to 0$ and $x_i(s_m)\to \ul x_i$. For $\vare>0$ small and $m$ large, we  get
$${{x_i'(s_m)}\over { \k_i(s_m)}}
\ge{{\ul x_i-\vare}\over {v_i\k_i(s_m)}}\left (\sum_{k=1}^{m_i} {{\al_{ik}(s_m)v_i}\over {1+\be_{ik}(s_m)(\ul x_i-\vare)}}+\sum_{j=1}^nd_{ij}(s_m) v_j\right )-{{\mu_i (s_m)}\over {\k_i(s_m)}}x_i(s_m)-x_i^2(s_m)$$
Since $\ul x_i\le v_iM_0$, taking limits $m\to\infty,\vare \to 0^+$, this leads to
$$0\ge \ul x_i\liminf_n \left [{1\over {v_i\k_i(s_m)}}\left (\sum_{k=1}^{m_i} {{\al_{ik}(s_m)v_i}\over {1+\be_{ik}(s_m)v_iM_0}}-\mu_i (s_m)v_i+\sum_{j=1}^nd_{ij}(s_m) v_j\right )-\ul x_i\right ],
$$
and therefore $\ul x_i/v_i\ge m_0$ for $m_0$ as in (5.16).\ter

\smal


Note that  when $\be_{ik}(t)\equiv 0$ in (5.8),  the  lower bound $m_0$ in (5.16) does not depend on $M_0$. 

In the case $n=1$, we get  a criterion which improves and  generalizes Theorem 3.2 in [9].

\proclaim{Corollary 5.1}. Consider the scalar equation (5.7), where $\al_k,\be_k,\mu,\k:[0,\infty)\to [0,\infty)$ are continuous and bounded, $\tau_k:[0,\infty)\to [0,\infty)$ are continuous, $1\le k\le m$, and $\k(t)>\k_0,A(t):=\sum_{k=1}^m \al_k(t)>A_0$ for some constants $\k_0,A_0>0$. If
$$ \inf_{t\ge 0}\left (\sum_{k=1}^m\al_k(t)\right )>\sup_{t\ge 0}\mu(t),
$$
then (5.7) is persistent (in $BC_0$) and all positive solutions are bounded. If in addition $\lim\limits_{t\to\infty} (t-\tau_{k}(t))=\infty$ for $ k=1,\dots,m$,
  (5.7) is permanent and positive solutions $x(t)=x(t;0,\var)\,  (\var \in BC_0)$ satisfy   
$$
m_0\le \liminf_{t\to\infty} x(t)\le \limsup_{t\to\infty} x(t)\le M_0,
$$
where 
$\dps 
M_0=\limsup_{t\to\infty}{1\over {\k(t)}} \left (\sum_{k=1}^{m}\al_{k}(t)-\mu(t)\right ),\
m_0=\liminf_{t\to\infty}{1\over {\k(t)}}\left (\sum_{k=1}^{m}{{\al_{k}(t)}\over {1+\be_{k}(t)M_0}}-\mu(t)\right ).
$

\bigskip

{\bf 5.3. Permanence for a competitive compartmental system with  infinite delay}

\med

In 1990, Aiello and Freedman [2] proposed and studied a model for a  single species with immature and mature stages and a discrete time-delay:
$$ \eqalign{
u_i'(t)&=\al u_m(t)-\gamma u_i(t)-\al e^{-\ga \tau}u_m(t-\tau)\cr
u_m'(t)&=\al e^{-\ga \tau} u_m(t-\tau)-\be u_m^2(t),\cr}
\eq(5.17)$$
where $\al,\be,\ga,\tau >0$. In (5.17), $u_i,u_m$ stand for the immature and mature populations, respectively, $\tau$ is the maturation time since birth, $\al$ is the birth rate for the species, and $\be, \ga$  the death rates for matures and immatures, respectively. Since the second equation is decoupled from the first, to describe the qualitative behavior of solutions to (5.17) it is sufficient to study the second equation of the system. This model has received great attention over the last decades,  a large number of generalizations has been derived, and many aspects of their dynamics analyzed, in several contexts. Besides the early works [2,3,13], see [4,5,28] and references therein for more results  on   related models.
A natural  extension is to introduce non-constant delay; in fact, it may even be infinite, leading to a second equation given by
$$u_m'(t)=\al \int_0^\infty f(s)e^{-\ga s} u_m(t-s)\, ds-\be u_m^2(t)$$ 
for some $f$ summable with $\int_0^\infty f(s)\, ds=1$. More recently, models with two or more stage-structured species  have been proposed. In [5], two species, structured  in matures and immatures and in competition, were considered. Disregarding the immature populations, this leads to a Lotka-Volterra type model of the form
$$ \eqalign{
u_1'(t)&=\al_1 \int_0^\infty f_1(s)e^{-\ga_1 s} u_1(t-s)\, ds-\be_1 u_1^2(t)-c_1u_1(t)u_2(t)\cr
u_2'(t)&=\al_2 \int_0^\infty f_2(s)e^{-\ga_2 s} u_2(t-s)\, ds-\be_2 u_2^2(t)-c_2u_1(t)u_2(t),\cr}
\eq(5.18)$$
where $\al_j,\be_j,\ga_j,c_j>0$ and the kernels $f_j:[0,\infty)\to [0,\infty)$ are continuous (or at least measurable) with $\int_0^\infty f_j(s)\, ds=1,\, j=1,2$. 
However, due to technical difficulties, in [4,5] the authors restricted their analyses to the case of kernels $f_1(s),f_2(s)$ with compact support -- of course,  if $[0,\tau]$ contains both the supports of $f_1$ and $f_2$, the framework is restricted to a system of differential equations with finite distributed delay on $[-\tau,0]$.  The theorem below was proven by Al-Omari and Gourley [5] for the case of a finite delay.


\proclaim{Theorem 5.5}. [5] Consider (5.18) with $\al_j,\be_j,\ga_j,c_j>0$, $f_j:[0,\infty)\to [0,\infty)$ continuous, $f_j(s)=0$ for $s\ge \tau$ for some $\tau \in (0,\infty)$, and
$\int_0^\tau f_j(s)\, ds=1,\, j=1,2$. Assume that
$$\eqalignno{
c_2\al_1 \int_0^\tau f_1(s)e^{-\ga_1 s}\, ds&< \be_1\al_2 \int_0^\tau f_2(s)e^{-\ga_2 s}\, ds,&(5.19)\cr
c_1\al_2 \int_0^\tau f_2(s)e^{-\ga_2 s}\, ds&< \be_2\al_1 \int_0^\tau f_1(s)e^{-\ga_1 s}\, ds.&(5.20)\cr}$$
Then $u^*=(u_1^*,u_2^*)$, where
$$\eqalign{
u_1^*&={1\over {\be_1\be_2-c_1c_2}}\left ( \be_1\al_2 \int_0^\tau f_2(s)e^{-\ga_2 s}\, ds-c_2\al_1 \int_0^\tau f_1(s)e^{-\ga_1 s}\, ds\right ),\cr
u_2^*&={1\over {\be_1\be_2-c_1c_2}}\left ( \be_2\al_1 \int_0^\tau f_1(s)e^{-\ga_1 s}\, ds-c_1\al_2 \int_0^\tau f_2(s)e^{-\ga_2 s}\, ds\right ),\cr}\eq(5.21)
$$
 is the unique positive equilibrium of (5.18), and $u^*$ is a global attractor of all  solutions to (5.18) with  initial conditions in $C_0:=\{ \var:[-\tau,0]\to\R_+^2\, |\,  \var$ is continuous and $\var (0)>0\}$.

Here, as an application of the techniques in Section 3, we prove that the above result is still valid for the case of $\tau=\infty$. 
We use the ideas and arguments in [5], inserting them  in the present framework, which enables us to deal with the infinite delay. 
Clearly (5.18) is not cooperative, and consequently Theorems 3.1 and 3.2 cannot be applied directly.  

Consider (5.18) with  kernels $f_1,f_2:[0,\infty)\to[0,\infty)$ continuous with $L^1$-norm equal to 1, and assume (5.19)-(5.20) with $\tau=\infty$.  Choose $\de >0$ such that
$$\eqalignno{
c_2\left (\de \be_1+\al_1 \int_0^\infty f_1(s)e^{-\ga_1 s}\, ds\right )&< \be_1\al_2 \int_0^\infty f_2(s)e^{-\ga_2 s}\, ds,&(5.22)\cr
c_1\left (\de \be_2+\al_2 \int_0^\infty f_2(s)e^{-\ga_2 s}\, ds\right )&< \be_2\al_1 \int_0^\infty f_1(s)e^{-\ga_1 s}\, ds.&(5.23)\cr}$$

Fix any solution $u(t)=u(t;\var)$ of (5.18) with $\var\in BC_0$. In $C_\ga^0$ or $UC_\ga$ with $0<\ga<\min(\ga_1,\ga_2)$, we write (5.18) in the form (3.1) with $n=2$, where $F=(F_1,F_2), G=(G_1,G_2)$ are linear  functions given by
$$F_i(\phi)=\al_i \int_0^\infty f_i(s)e^{-\ga_i s} \phi_i(-s)\, ds\q {\rm for}\q \phi=(\phi_1,\phi_2)\in BC, i=1,2,$$
and $G(\phi)= G(\phi_1(0),\phi_2(0))$, with
$$G_1(u)=\be_1u_1+c_1u_2,\  G_2(u)=\be_2u_2+c_2u_1\q {\rm for}\q u=(u_1,u_2)\in \R^2.$$
Since $F(\phi)\ge 0$ for $\phi\in BC_+$, from Lemma 3.1 the solutions $u(t)$ of (5.18) with initial conditions in $BC_0$ are  non-negative in their maximal interval of existence. For such solutions, we derive
$$ \eqalign{
u_1'(t)&\le\al_1 \int_0^\infty f_1(s)e^{-\ga_1 s} u_1(t-s)\, ds-\be_1 u_1^2(t)\cr
u_2'(t)&\le\al_2 \int_0^\infty f_2(s)e^{-\ga_2 s} u_2(t-s)\, ds-\be_2 u_2^2(t),\cr}
$$
and we now compare with solutions to the auxiliary cooperative system
$$ \eqalign{
u_1'(t)&=\al_1 \int_0^\infty f_1(s)e^{-\ga_1 s} u_1(t-s)\, ds-\be_1 u_1^2(t)\cr
u_2'(t)&=\al_2 \int_0^\infty f_2(s)e^{-\ga_2 s} u_2(t-s)\, ds-\be_2 u_2^2(t).\cr}
\eq(5.24_{1,u})$$
System $(5.24_{1,u})$ satisfies (A1)-(A6), so in particular it is dissipative in $BC_0$. Thus, $u(t)$ satisfies
$$\limsup_{t\to\infty} u_i(t)\le u_i^{1,u},\, i=1,2,$$ where
$u^{1,u}=(u_1^{1,u},u_2^{1,u})$ is the unique positive equilibrium of $(5.24_{1,u})$, with coordinates given by
$$u_i^{1,u}={{\al_i}\over {\be_i}}\int_0^\infty f_i(s)e^{-\ga_i s}\, ds,\q i=1,2.$$
For $t>0$ sufficiently large, we get
$$ \eqalign{
u_1'(t)&\ge\al_1 \int_0^\infty f_1(s)e^{-\ga_1 s} u_1(t-s)\, ds-\be_1 u_1^2(t)-c_1u_1(t)(u_2^{1,u}+\de)\cr
u_2'(t)&\ge\al_2 \int_0^\infty f_2(s)e^{-\ga_2 s} u_2(t-s)\, ds-\be_2 u_2^2(t)-c_2u_2(t)(u_1^{1,u}+\de),\cr}
$$
and compare $u(t)$ with solutions to the new auxiliary  system
$$ \eqalign{
u_1'(t)&=\al_1 \int_0^\infty f_1(s)e^{-\ga_1 s} u_1(t-s)\, ds-\be_1 u_1^2(t)-c_1u_1(t)(u_2^{1,u}+\de)\cr
u_2'(t)&=\al_2 \int_0^\infty f_2(s)e^{-\ga_2 s} u_2(t-s)\, ds-\be_2 u_2^2(t)-c_2u_2(t)(u_1^{1,u}+\de).\cr}
\eq(5.24_{1,l})$$
System $(5.24_{1,l})$ has the form (3.1) and satisfies (A1)-(A4) and (A6). Also, for $(5.24_{1,l})$ the matrix $B$   in (A5) reads as $B=\pmatrix{c_1(u_2^{1,u}+\de)&0\cr
0&c_2(u_1^{1,u}+\de)\cr}$.  Formulae (5.22), (5.23) yield
$$\eqalign {&F(u)-Bu\cr
&=\Big (\al_1 u_1 \int_0^\infty f_1(s)e^{-\ga_1 s}\, ds-c_1(u_2^{1,u}+\de)u_1,\al_2 u_2 \int_0^\infty f_2(s)e^{-\ga_2 s}\, ds)-c_2(u_1^{1,u}+\de)u_2\Big )>0,\cr}$$ 
for any $u=(u_1,u_2)>0$. In particular,
(A5) holds. From Theorem 3.2, we now obtain 
$$\liminf_{t\to\infty} u_i(t)\ge u_i^{1,l},\, i=1,2,$$ where
$u^{1,l}=(u_1^{1,l},u_2^{1,l})$ is the unique positive equilibrium of $(5.24_{1,l})$, given by
$$u_1^{1,l}={{\al_1  \int_0^\infty f_1(s)e^{-\ga_1 s}\, ds-c_1(u_2^{1,u}+\de)}\over {\be_1}},
\ u_2^{1,l}={{\al_2  \int_0^\infty f_2(s)e^{-\ga_2 s}\, ds)-c_2(u_1^{1,u}+\de)}\over {\be_2}}.$$

Next, we  observe that, for $t>0$ sufficiently large, the solution $u(t)$ satisfies
$$ \eqalign{
u_1'(t)&\le\al_1 \int_0^\infty f_1(s)e^{-\ga_1 s} u_1(t-s)\, ds-\be_1 u_1^2(t)-c_1u_1(t)(u_2^{1,l}-\de)\cr
u_2'(t)&\le\al_2 \int_0^\infty f_2(s)e^{-\ga_2 s} u_2(t-s)\, ds-\be_2 u_2^2(t)-c_2u_2(t)(u_1^{1,l}-\de),\cr}
$$
and compare it with solutions of the second-upper auxiliary system
$$ \eqalign{
u_1'(t)&=\al_1 \int_0^\infty f_1(s)e^{-\ga_1 s} u_1(t-s)\, ds-\be_1 u_1^2(t)-c_1u_1(t)(u_2^{1,l}-\de)\cr
u_2'(t)&=\al_2 \int_0^\infty f_2(s)e^{-\ga_2 s} u_2(t-s)\, ds-\be_2 u_2^2(t)-c_2u_2(t)(u_1^{1,l}-\de)\cr},
\eq(5.24_{2,u})$$
which satisfies (A1)-(A6).
Proceeding in this way, in [5] the arguments above were iterated and  auxiliary cooperative systems $(5.24_{n,u}),(5.24_{n,l})$  with positive equilibria $u^{n,u}, u^{n,l}$, respectively, constructed,  providing upper and lower bounds for
$u(t)$:
$$ u_i^{n,l}\le \liminf_{t\to\infty} u_i(t)\le \limsup_{t\to\infty} u_i(t)\le u_i^{n,u},\q i=1,2,\ n\in \N.$$
Moreover,  explicit recursive formulae
can be obtained  for the sequences $ u_i^{n,l}, u_i^{n,u}$, showing  that $ u_i^{n,l}$  increases and $u_i^{n,u}$ decreases, $i=1,2$, with
$\lim_n u^{n,l}=\lim_n u^{n,u}=u^*.$
See [5] for details.  Thus, 

\proclaim{Theorem 5.6}. For (5.18), Theorem 5.5 is valid with $\tau=+\infty$ and $C_0$ replaced by $BC_0$.

\med

{\bf Acknowledgement}: 
Work  supported by Funda\c c\~ao para a Ci\^encia e a Tecnologia, under 
 PEst-OE/MAT/UI0209/2013.

\bigskip

{\bf References}

\baselineskip=13.5pt

\med


\item{1.}   Ahmad, S., and Lazer, A.C. (2006). Average growth and total permanence in a competitive
Lotka-Volterra system. Ann. Mat. Pura Appl.  185, S47--S67.

\item{2.}  Aiello, W.G., and  Freedman, H.I. (1990). A time-delay model of single species growth with stage structure. Math. Biosci. 101, 139--153.

\item{3.} Aiello, W.G.,   Freedman, H.I., and  Wu, J. (1992). Analysis of a model representing stage-structured population growth with state-dependent time delay. SIAM J.~Appl.~Math.~52, 855--869.

\item{4.}  Al-Omari, J.F.M., and Al-Omari, S.K.Q. (2011). Global stability in a structured population model with distributed maturation delay and harvesting. Nonlinear Anal. RWA 12, 1485--1499.

\item{5.} Al-Omari, J.F.M., and  Gourley, S.A. (2003). Stability and traveling fronts in Lotka-Volterra competition models with stage structure. SIAM J. Appl. Math. 63, 2063--2086.

\item{6.} Arino, J., Wang,  L.,  and  Wolkowicz, G.S.K. (2006), An alternative formulation for a delayed logistic equation.  {J. Theor. Biol.} { 241}, 109--119.
   
   \item{7.}
     Bastinec, J.,  Berezansky, L.,  Diblik, J., and   Smarda, Z.  (2014). On a delay population model with a quadratic nonlinearity without positive steady state. Appl.~Math.~Comput.~227, 622--629.

%
%
%

\item{8.} 
 Faria, T. (2010). Stability and extinction for Lotka-Volterra systems with infinite delay. J. Dynam. Differential Equations 22, 299--324.


\item{9.} Faria, T. (2014).  A note on permanence of nonautonomous cooperative scalar population models  with delays.  {Appl. Math. Comput.} 240, 82--90.

\item{10.} Faria, T. (2014). Global dynamics for  Lotka-Volterra Systems with  infinite delay and patch structure.  {Appl. Math. Comput.} 245, 575--590.

\item{11.}
 Faria, T.,  and  Muroya, Y. Global attractivity and extinction for  Lotka-Volterra systems with infinite delay and feedback controls. Proc. Roy. Soc. Edinburgh Sect. A, {\it to appear}. (arXiv:1307.7039v1 [math.CA]).

 
\item{12.}  
 Faria, T., and Oliveira, J.J.  (2011). General criteria for asymptotic and exponential stabilities of neural network models with unbounded delay. Appl. Math. Comput.  217, 9646--9658.

 \item{13.}  Freedman, H.I., and  Wu, J.H. (1991). Persistence and global asymptotic stability of single species dispersal models with stage structure.
Quart. Appl. Math. 49, 35--371. 


\item{14.} Haddock, J.R., and Hornor, W. (1988). Precompactness and convergence in norm of positive orbits in a certain
fading memory space.  Funkcial. Ekvac. 31, 349--361.

\item{15.}  Haddock, J.R.,  Nkashama, M.N.,  Wu, J.H. (1989).
    Asymptotic constancy for linear neutral Volterra integrodifferential
    equations, T\^{o}hoku Math. J. 41, 689-710. 

\item{16.} Hale, J.K. (1974). Functional differential equations with infinite delay. J. Math Anal. Appl. 48, 276--283.

\item{17.}
Hale, J.K., and   Kato, J. (1978). Phase space for retarded equations with infinite delay, Funkcial.
Ekvac. 21, 11--41.



\item{18.}
 Hino, Y.,   Murakami, S., and   Naito, T., (1993). Functional Differential
Equations with Infinite Delay, Sprin\-ger-Verlag, New-York.


\item{19.}  Kuang, Y. (1993).  Delay Differential Equations
with Applications in Population Dynamics, Academic Press, New York.

\item{20.}  Kuang, Y. (1995).  Global stability in delay differential systems without dominating instantaneous
negative feedbacks.  J. Differential Equations 119, 503--532.

\item{21.}  Kuang, Y., and  Smith, H.L. (1993). Global stability for infinite delay Lotka-Volterra type systems.
 J. Differential Equations {103}, 221--246.


\item{22.}  Mallet-Paret, J. and  Sell, G.R. (1996).
Systems of differential delay equations:
Floquet multipliers and discrete Lyapunov functions.
J. Differential Equations 125, 385--440.

\item{23.} Mallet-Paret, J. and  Sell, G.R. (1996). The Poincar\'e-Bendixson theorem for monotone
cyclic feedback systems with delay. J. Differential Equations 125, 441--489.

 \item{24.}  Meng, X., and  Chen, L. (2008). Periodic solution and almost periodic solution for a nonautonomous LotkaÐVolterra dispersal system with infinite delay,  J.~Math.~Anal.~Appl. 339, 125--145.


\item{25.}  Miller, R.K. (1966). Asymptotic behavior of solutions of nonlinear Volterra equations. Bull. Amer. Math. Soc. 72, 153Ð156.

\item{26.} Miller, R.K., and Sell, G.R. (1968).
Existence, uniqueness and continuity of solutions of integral equations.
Ann. Mat. Pura Appl. 80, 135--152.

\item{27.}  Murakami, S., and  Naito,  T. (1989). Fading memory spaces and stability properties for functional differential equations with infinite delay. Funkcial. Ekval. 32, 91--105.


\item{28.} Saito, Y., and Takeuchi, Y. (2003). A time-delay for predator-prey growth with stage structure. Cann. Appl. Math. Quart. 13, 293--302.

\item{29.}  Smith, H.L. (1986). Invariant curves for mappings. SIAM J. Math. Anal. 15, 1053--1067.

\item{30.}
Smith, H.L. (1995). Monotone Dynamical Systems. An
Introduction to the Theory of Competitive and Cooperative Systems, Amer. Math.  Soc., Providence, RI. 
 
 \item{31.}
 Smith, H.L. and  Thieme,  H.R. (2011). Dynamical Systems and Population Persistence,  
Amer. Math. Soc., Providence, RI.




 \item {32.}   Teng, Z., and Chen, L. (2001). Global stability of periodic Lotka-Volterra systems with delays. Nonlinear Anal. 45, 1081--1095.
   
 \item {33.}   Teng, Z., and   Rehim,  M. (2006). Persistence in nonautonomous predator-prey systems with infinite delay. J. Comput. Appl. Math.  197, 302--321.

 
 \item{34.} Wu,  J. (1992). Global dynamics of strongly monotone retarded equations with infinite delay. J. Integral Equations Appl. 4, 273--307.

%
%
%

\item{35.}
Zhao, X.-Q., and Jing, Z.-J. (1996). Global asymptotic behavior in some cooperative systems of functional differential equations. Cann. Appl. Math. Quart. 4, 421-444. 
\end

\end